\documentclass{amsart}
\usepackage{latexsym,amsmath,amssymb}

\theoremstyle{definition}

\theoremstyle{remark}

\numberwithin{equation}{section}

\begin{document}
\title[ON THE SPACE $\ell^{p}(\beta)$]{\textsc{\textbf{ON THE SPACE $%
\ell^{p}(\beta)$ WITH A WEIGHTED CAUCHY PRODUCT }}}
\author{\textsc{Y. Estaremi, M. R. Jabbarzadeh and M. T. Karaev,}}
\address{\textsc{y. estaremi$^{1}$ m. r. jabbarzadeh$^{2}$, and m. t. karaev$%
^{3}$}} \email{estaremi@gmail.com}
\email{mjabbar@tabrizu.ac.ir}\email{mubariztapdigoglu@sdu.edu.tr}

\address{$^{1}$  Suleyman Demirel University, Isparta Vocational School,
32260 Isparta, Turkey\\
and\\
Institute of Mathematics and Mechanics National Academy of Sciences of
Azerbaijan, F.Agayev 9, Baku, Azerbaijan\\
$^{2,3}$ Faculty of Mathematical Sciences, University of Tabriz, P. O. Box:
5166615648, Tabriz, Iran}
\thanks{}
\thanks{}
\date{}
\dedicatory{}

\begin{abstract}
In this paper we consider a weighted Cauchy product $\diamond$ on $%
\ell^{p}(\beta)$ and then we characterized some Banach algebra structures
for $\ell^{p}(\beta)$.
\end{abstract}

\maketitle

\commby{}

%%% ----------------------------------------------------------------------

\section{\textsc{Introduction}}

\vspace*{0.3cm}Let $\{\beta_{n}\}^{\infty}_{n=0}$ be a sequence of positive
numbers with $\beta(0)=1$. For $1\leq p <\infty$ we consider the space of
sequences $f=\{\hat{f}(n)\}$ with $\|f\|_{\beta}^p=\Sigma^{\infty}_{n=0}|
\hat{f}(n)|^{p}\beta(n)^{p}<\infty$. We shall use the formal notation $f(z)=
\Sigma^{\infty}_{n=0}\hat{f}(n)z^{n}$ whether or not the series converges
for any complex values of $z$. Let $\ell^{p}(\beta)$ denote the space of
such formal power series.

\vspace*{0.3cm} After the first investigations of Nikol'skii(see
in \cite{nik} and A. L. Shields \cite{sh}, several publications
have appeared on weighted shift operators.,see,for instance,
Lamber\cite{la, lam1, lam2, lam3},Domar\cite{do}(and references
therein),Yadav\cite{yada1, yada2} and Karaev\cite{kar1}. In this
paper we consider a weighted Cauchy product $\diamond $ on $\ell
^{p}(\beta )$ and then, by another conditions or calculation
methods , we characterized some classic Banach algebra structures
for $\ell ^{p}(\beta )$ with $1\leq p<\infty $.

\vspace*{0.3cm} In the next section, we define a weighted Cauchy product $%
\diamond$, under certain conditions, on $\ell^{p}(\beta)$ and then we show
that the Banach space $\ell^{p}(\beta)$ with weighted Cauchy product $%
\diamond$ is a Banach algebra. Then we determine invertible elements and
maximal ideal space of $(\ell^{p}(\beta),\diamond)$. Also, we give a
sufficient condition for the $\diamond$-multiplication operator $%
M_{\diamond,z}$ acting on $\ell^{p}(\beta)$ to be unicellular.  This result,
as usual, leads to a description of closed ideals of the algebra $%
\ell^{p}(\beta)$ and cyclic vectors of the $\diamond$-multiplication
operator $M_{\diamond,z}$.

\section{ \textsc{Some Banach Algebra Structures For $\ell^{p}(\protect\beta)
$}}

Let $\{\delta_{n}\}^{\infty}_{n=0}$ be  a sequence of positive numbers with $%
\delta_{0}=1.$ Let $1<p<\infty$ and let $q$ be the conjugate exponent to $p$%
. For each $k, M\in\mathbb{N}\cup\{0\}$, take
\begin{equation}
C_o:=\sup_{n\geq 0}\sum^{n}_{k=0}\left(\frac{\delta_{n}\beta(n)} {%
\delta_{k}\delta_{n-k}\beta(k)\beta(n-k)}\right)^{q},
\end{equation}
and
\begin{equation}
b_{M, k}:=\sup_{n\geq M+1}\frac{\delta_{n+k}\beta(n+k)}{\delta_{n}\delta_{k}%
\beta(n)\beta(k)} \ .
\end{equation}
%Throughout this paper we assume that $1<p<\infty$,  $C_o<\infty$
%and $\lim b_{M, k}=0$ when $M\rightarrow\infty$.
Given arbitrary two functions $f(z)=\sum^{\infty}_{n=0}\hat{f}(n)z^{n}$  and
$g(z)=\sum^{\infty}_{n=0}\hat{g}(n)z^{n}$ of the space $\ell^{p}(\beta)$,
define the following weighted Cauchy product series
\begin{eqnarray}
\ f\diamond \ g=\sum^{\infty}_{n=0}\sum^{\infty}_{m=0}\frac{\delta_{m+n}} {%
\delta_{n}\delta_{m}}\hat{f}(n)\hat{g}(m)z^{m+n}.
\end{eqnarray}

Note that if we set $\delta _{n}\equiv 1$, the weighted Cauchy product $%
\diamond $ will be coincide to the ordinary Cauchy product. Also
weighted Cauchy product $\diamond $, can be regarded as a discrete
case od Duhamel product (see \cite{wig}). For more study of
Duhamel product on various function spaces and its applications,
see \cite{kar1, kar2, kar3} and the references therein.

\vspace*{0.3cm} Let $\ell ^{0}(\beta )$ be the set of all formal power
series. For each $f\in \ell ^{p}(\beta )$, let $M_{\diamond ,f}:\ell
^{p}(\beta )\rightarrow \ell ^{0}(\beta )$ defined by $M_{\diamond
,f}(g)=f\diamond g$ be its corresponding $\diamond $-multiplication linear
operator. It is easy to see that $M_{\diamond ,z}(f)=\sum_{n=0}^{\infty }%
\frac{\delta _{n+1}}{\delta _{n}\delta _{1}}\hat{f}(n)z^{n+1}$ and $%
M_{\diamond ,z}^{N}(f)=\frac{\delta _{N}}{\delta _{1}^{N}}z^{N}\diamond f$,
for all $N\in \mathbb{N}\cup \{0\}$ and $f\in \ell ^{p}(\beta )$. Take $%
C=\sup_{n\geq 0}\frac{\beta (n+N)\delta _{n+N}}{\beta (n)\delta _{n}\delta
_{1}^{N}}$. If for each $k,M\in \mathbb{N}\cup \{0\}$, $b_{M,k}<\infty $,
then $C<\infty $ and in this case $\Vert M_{\diamond ,z}^{N}\Vert =C$ (see
\cite{sh}). Note that $M_{\diamond ,z}$ on $(\ell ^{p}(\beta ),\diamond )$
is equal to $M_{z}$ on $\ell ^{p}(\widetilde{\beta })$ with ordinary Cauchy
product, where $\widetilde{\beta }(n):=\frac{\delta _{n+1}\beta (n)}{\delta
_{n}\delta _{1}}$. However, $\ell ^{p}(\beta )$ is not unitarily equivalent
to $\ell ^{p}(\widetilde{\beta })$, because in general the sequence $\{\frac{%
\widetilde{\beta }(n)}{\beta (n)}\}$ is not necessarily constant .

%. It follows that
% $$\|M^N_{\diamond,z}f\|^{p}_{\beta}=\sum^{\infty}_{n=0}(\frac{\delta_{n+N}}
%{\delta_{n}\delta^N_{1}})^{p}|\hat{f}(n)|^{p}\beta(n+N)^{p}$$$$=\sum^{\infty}_{n=0}|\hat{f}(n)|
% ^{p}\beta(n)^{p}(\frac{\delta_{n+N}\beta(n+N)}{\delta_{n}\beta(n)\delta^N_{1}})^{p}
% \leq C^{p}\|f\|^{p}_{\beta} .$$
%Hence $\|M^N_{\diamond,z}\|\leq C$. On the other hand if we put
%$f_{n}(z)=z^{n}$, then $\|f_{n}\|_{\beta}=\beta(n)$ and
%$M^N_{\diamond,z}(f_{n})
%=\frac{\delta_{n+N}}{\delta_{n}\delta^N_{1}}z^{n+N}$. Therefore we
%have
% $$\frac{\beta(n+N)\delta_{n+N}}{\delta_{n}\delta^N_{1}}=
% \|M^N_{\diamond,z}(f_{n})\|_{\beta}\leq\|M^N_{\diamond,z}
% \| \|f_{n}\|_{\beta}=\|M^N_{\diamond,z}\| \beta(n).$$
% This implies that $C\leq\|M^N_{\diamond,z}\|$ and so $\|M^N_{\diamond,z}\|=C$ (see \cite{sh}).

\vspace*{0.3cm} %{\bf Theorem 2.1.} $(\ell^{p}(\beta),\diamond)$
%is a unital commutative Banach algebra.

%\vspace*{0.3cm} {\bf Proof.}
Now, if (2.1) is holds, then by using H\"{o}lder inequality we get that $%
\Vert M_{\diamond ,f}\Vert \leq {C_{o}}^{1/q}\Vert f\Vert _{\beta }$, for
each $f\in \ell ^{p}(\beta )$. Also it is easy to see that the constant
function $f=1$ is a unity for $(\ell ^{p}(\beta ),\diamond )$. Therefore $%
(\ell ^{p}(\beta ),\diamond )$ is a unital commutative Banach algebra . In
the following proposition we consider the case $p=1$.

%it is sufficient prove that $M_{\diamond, f}$ is a bounded
%$f,g\in\ell^{p}(\beta)$. Using (2.3), it is easy to see that
%\begin{equation}
%\widehat{({f\diamond} g)}(n)=\sum^{n}_{k=0}\frac{\delta_{n}}
%{\delta_{k}\delta_{n-k}}\hat{f}(k)\hat{g}(n-k).
%5\end{equation}

%By using H\"{o}lder inequality and (2.1) we have
%$$
%\|M_{\diamond,
%f}(g)\|_{\beta}^p=\sum^{\infty}_{n=0}|\widehat{(f\diamond
%g)}(n)|^{p}\beta(n)^{p}$$$$=\sum^{\infty}_{n=0}\left|\sum^{n}_{k=0}\frac{\delta_{n}}{\delta_{k}
%\delta_{n-k}}\hat{f}(k)\hat{g}(n-k)\right|^{p}\beta(n)^{p}$$$$
%\leq\sum^{\infty}_{n=0}\left(\sum^{n}_{k=0}\frac{\delta_{n}\beta(n)}
%{\delta_{k}\delta_{n-k}
%\beta(k)\beta(n-k)}|\hat{f}(k)||\hat{g}(n-k)|\beta(k)\beta(n-k)\right)^{p}
%$$
%$$\leq\sum^{\infty}_{n=0}\left(\sum^{n}_{k=0}\left(|\hat{f}(k)|\beta(k)|
%\hat{g}(n-k)|\beta(n-k)\right)^{p}\right)^{\frac{p}{p}}
%\left(\sum^{n}_{k=0}\left(\frac{\delta_{n}\beta(n)}{\delta_{k}
%\delta_{n-k}\beta(k)\beta(n-k)}\right)^{q}\right)^{\frac{p}{q}}
%$$$$\leq C_{o}^{\frac{p}{q}}\sum^{\infty}_{n=0}\sum^{n}_{k=0}|\hat{f}(k)|^{p}\beta(k)^{p}|
%\hat{g}(n-k)|^{p}\beta(n-k)^{p}$$$$=
%C_{o}^{\frac{p}{q}}\left(\sum^{\infty}_{n=0}|\hat{f}(n)|^{p}\beta(n)^{p}\right)\left(\sum^{\infty}_{n=0}
%|\hat{g}(n)|^{p}\beta(n)^{p}\right)=C_{o}^{\frac{p}{q}}\|f\|_\beta^p \ \|g\|_\beta^p.$$

%Consequently, we get that
%$$\|M_{\diamond, f}(g)\|_{\beta}=\|f\diamond
%g\|_{\beta}\leq{C_o}^{\frac{1}{q}}\|f\|_{\beta}  \|g\|_{\beta},$$
%and so $\|M_{\diamond, f}\|\leq
%{C_o}^{\frac{1}{q}}\|f\|_{\beta}$. \ $\Box$

\vspace*{0.3cm}\textbf{Proposition 2.1.} Let for some $N\in\mathbb{N}$,
\begin{equation}
\sum_{n,m\geq N+1}\frac{\delta_{n+m}\beta(n+m)}{\delta_{n}\delta_{m}\beta(n)%
\beta(m)}<\infty .
\end{equation}
Then $(\ell^{1}(\beta),\diamond)$ is a unital commutative Banach algebra.

\vspace*{0.3cm} \textbf{Proof.} Let $f\in\ell^{1}(\beta)$. It is sufficient
prove that $M_{\diamond, f}$ is a bounded operator on $\ell^{1}(\beta)$. To
see this let $g\in\ell^{1}(\beta)$. Then we have

\begin{equation*}
f\diamond g=\sum^{\infty}_{n=0}\sum^{\infty}_{m=0}\frac{\delta_{m+n}} {%
\delta_{n}\delta_{m}}\hat{f}(n)\hat{g}(m)z^{m+n}=\hat{f}(0)g+\frac{\hat{f}(1)%
}{\delta_{1}}\sum^{\infty}_{m=0}\frac{\delta_{m+1}} {\delta_{m}}\hat{g}%
(m)z^{m+1}\ \ + \ \ \cdots
\end{equation*}

\begin{equation*}
+\frac{\hat{f}(N)}{\delta_{N}}\sum^{\infty}_{m=0}\frac{\delta_{m+N}} {%
\delta_{m}}\hat{g}(m)z^{m+N}+\hat{g}(0)\sum^{\infty}_{n=N+1}\hat{g}(n)z^{n}+%
\frac{\hat{g}(1)}{\delta_{1}}\sum^{\infty}_{n=N+1}\frac{\delta_{n+1}} {%
\delta_{n}}\hat{f}(n)z^{n+1}+\cdots
\end{equation*}
\begin{equation*}
+\frac{\hat{g}(N)}{\delta_{N}}\sum^{\infty}_{n=0}\frac{\delta_{n+N}} {%
\delta_{n}}\hat{f}(n)z^{m+N}+\sum^{\infty}_{n=N+1}\sum^{\infty}_{m=N+1}\frac{%
\delta_{m+n}} {\delta_{n}\delta_{m}}\hat{f}(n)\hat{g}(m)z^{m+n}.
\end{equation*}

Thus we can write

\begin{equation*}
f\diamond g=\sum^{\infty}_{n=0}\sum^{\infty}_{m=0}\frac{\delta_{m+n}} {%
\delta_{n}\delta_{m}}\hat{f}(n)\hat{g}(m)z^{m+n}=\hat{f}(0)g+\hat{f}%
(1)M_{z,\diamond}(g)+  \cdots
\end{equation*}
\begin{equation*}
+\frac{\hat{f}(N)\delta^{N}_{1}}{\delta_{N}}M^{N}_{z,\diamond}(g)  +\hat{g}%
(0)R_{N+1}(f)+\hat{g}(1)R_{N+1}(M_{z,\diamond}(f))+  \cdots
\end{equation*}%
\begin{equation*}
+\frac{\hat{f}(N)\delta^{N}_{1}}{\delta_{N}}R_{N+1}(M^{N}_{z,\diamond}(g))+%
\sum^{\infty}_{n=N+1}\sum^{\infty}_{m=N+1}\frac{\delta_{m+n}} {%
\delta_{n}\delta_{m}}\hat{f}(n)\hat{g}(m)z^{m+n},
\end{equation*}

where $R_{N}(g):=\sum^{\infty}_{n=N}\hat{g}(n)z^{n}$.  It follows that

\begin{equation*}
\|f\diamond g\|_{\beta}\leq\left(2\sum^{N}_{k=0}\frac{\delta^k_{1}}{%
\beta(k)\delta_{k}}\|M^k_{\diamond,z}\|  +\sum^{\infty}_{n,m=N+1}\frac{%
\delta_{m+n}\beta(n+m)} {\delta_{n}\delta_{m}\beta(n)\delta(m)}
\right)\|f\|_{\beta}\|g\|_{\beta},
\end{equation*}
and so $M_{\diamond, f}$ is bounded.  \ $\Box$

\vspace*{0.3cm} From now on we assume that for $p=1$, (2.4) and for $%
1<p<\infty$, (2.1) and (2.2) are holds.

\vspace*{0.3cm} \textbf{Lemma 2.2.} Let $f\in\ell^{p}(\beta)$ and let $%
\lambda$ be a nonzero complex number. If $\hat{f}(0)=0$ then $\lambda
I-M_{\diamond, f}$ has closed range, where $I$ is the identity operator.

\vspace*{0.3cm} \textbf{Proof.} Let $\hat{f}(0)=0$. To show that $\lambda
I-M_{\diamond, f}$ has closed range, we only need to prove the $\diamond$%
-multiplication operator $M_{\diamond, f}$ is compact on $\ell^{p}(\beta)$.
For $M\in\mathbb{N}$ define $K_{M}$ on $\ell^{p}(\beta)$ by
\begin{equation*}
K_{M}(g)=\sum^{M}_{m=0}\sum^{\infty}_{n=0}\frac{\delta_{n+m}} {%
\delta_{n}\delta_{m}}\hat{f}(n)\hat{g}(m)z^{n+m}.
\end{equation*}
Since
\begin{equation*}
\sum^{M}_{m=0}\sum^{\infty}_{n=M+1}\frac{\delta_{n+m}} {\delta_{n}\delta_{m}}%
\hat{f}(n)\hat{g}(m)z^{n+m}=\sum^{M}_{m=0}\frac{\hat{g}(m)\delta^m_{1}}{%
\delta_{m}} \left(\sum^{\infty}_{n=M+1}\frac{\delta_{n+m}} {%
\delta_{n}\delta^m_{1}}\hat{f}(n)z^{n+m}\right)
\end{equation*}
\begin{equation*}
= \sum^{M}_{m=0} \frac{\hat{g}(m)\delta^m_{1}}{\delta_{m}}%
R_{M+1}(M^m_{\diamond,z}(f)) ,
\end{equation*}
we have
\begin{equation*}
K_{M}(g):=\sum^{M}_{m=0}\sum^{M}_{n=0}\frac{\delta_{m+n}}{%
\delta_{n}\delta_{m}} \hat{f}(n)\hat{g}(m)z^{n+m}+\sum^{M}_{m=0} \frac{\hat{g%
}(m)\delta^m_{1}}{\delta_{m}}R_{M+1}(M^m_{\diamond,z}(f)) ,
\end{equation*}
and so $K_{M}$ is a bounded and finite-rank operator on $\ell^{p}(\beta)$.
Also, it is easy to verify that
\begin{equation*}
\sum^{\infty}_{n=M+1}\sum^{\infty}_{m=M+1}\frac{\delta_{n+m}}{\delta_{n}
\delta_{m}}\hat{g}(n)\hat{f}(m)z^{n+m}=\sum^{\infty}_{n=2M+2}(%
\sum^{n-M-1}_{k=M+1}\frac{\delta_{n}}{\delta_{k} \delta_{n-k}}\hat{g}(k)\hat{%
f}(n-k))z^{n}.
\end{equation*}
Then for $1<p<\infty$ we have
\begin{equation*}
\|M_{\diamond, f}(g)-K_{M}(g)\|_{\beta}=\|f\diamond g-K_{M}(g)\|_{\beta}
\end{equation*}
\begin{equation*}
=\left
\|\sum^{\infty}_{n=M+1}\sum^{\infty}_{m=M+1}\frac{\delta_{n+m}}{%
\delta_{n} \delta_{m}}\hat{g}(n)\hat{f}(m)z^{n+m}+
\sum^{M}_{n=0}\sum^{\infty}_{m=M+1}\frac{\delta_{n+m}} {\delta_{n}\delta_{m}}%
\hat{f}(n)\hat{g}(m)z^{n+m}\right\|_{\beta}
\end{equation*}
\begin{equation*}
\leq\left\|\sum^{\infty}_{n=2M+2}(\sum^{n-M-1}_{k=M+1}\frac{\delta_{n}}{%
\delta_{k} \delta_{n-k}}\hat{g}(k)\hat{f}(n-k))z^{n}\right\|_{\beta}
\end{equation*}
\begin{equation*}
+\frac{|\hat{f}(1)|}{\delta_{1}}\left(\sum^{\infty}_{m=M+1}\frac{%
\delta^{p}_{m+1}} {\delta^{p}_{m}}|\hat{g}(m)|^{p}\beta(m+1)^{p}\right)^{%
\frac{1}{p}}
\end{equation*}%
\begin{equation*}
+\cdots+ \frac{|\hat{f}(M)|}{\delta_{M}}\left(\sum^{\infty}_{m=M+1}\frac{%
\delta^{p}_{m+M}} {\delta^{p}_{m}}|\hat{g}(m)|^{p}\beta(m+M)^{p}\right)^{%
\frac{1}{p}} .
\end{equation*}

Since for every $0\leq k\leq M$,

\begin{equation*}
\frac{|\hat{f}(k)|}{\delta_{k}}\left(\sum^{\infty}_{m=M+1}\frac{%
\delta^{p}_{m+k}} {\delta^{p}_{m}}|\hat{g}(m)|^{p}\beta(m+k)^{p}\right)^{%
\frac{1}{p}}
\end{equation*}%
\begin{equation*}
= |\hat{f}(k)|\beta(k)\left(\sum^{\infty}_{m=M+1}\frac{\delta^{p}_{m+k}%
\beta(m+k)^{p}} {\delta^{p}_{m}\delta^{p}_{k}\beta(k)^{p}\beta(m)^{p}}|\hat{g%
}(m)|^{p}\beta(m)^{p}\right)^{\frac{1}{p}}
\end{equation*}
\begin{equation*}
\leq\|f\|_\beta\left(\sup_{m\geq{M+1}}\frac{\delta_{m+k}\beta(m+k)} {%
\delta_{m}\delta_{k}\beta(m)\beta(k)}\right)\left(\sum_{m=1}^{\infty}|\hat
g(m)|^p\beta(m)^p\right)^{\frac{1}{p}} \leq
\|f\|_{\beta}\|g\|_{\beta}b_{M,k},
\end{equation*}
then we get that

\begin{equation*}
\|M_{\diamond, f}(g)-K_{M}(g)\|_{\beta}
\end{equation*}%
\begin{equation*}
\leq \left(\sum^{\infty}_{n=2M+2}\left(\sum^{n-M-1}_{k=M+1}\frac{%
\delta_{n}\beta(n)}{\delta_{k} \delta_{n-k}\beta(k)\beta(n-k)}|\hat{g}%
(k)|\beta(k)|\hat{f}(n-k)|\beta(n-k)\right)^{p}\right)^{\frac{1}{p}}
\end{equation*}
\begin{equation*}
+\|f\|_{\beta}\|g\|_{\beta}(b_{M,1}+b_{M,2}+\cdots+b_{M,M})
\end{equation*}
\begin{equation*}
\overset{\mathrm{H\ddot{o}lder}}{\leq}\left(\sum^{\infty}_{n=2M+2}\left(%
\sum^{n-M-1}_{k=M+1}|\hat{f}(k)|^{p}\beta(k)^{p}| \hat{g}(n-k)|^{p}%
\beta(n-k)^{p}\right)\right.
\end{equation*}%
\begin{equation*}
\times\left
.\left(\sum^{n-M-1}_{k=M+1}\left(\frac{\delta_{n}\beta(n)}{%
\delta_{k}\delta_{n-k} \beta(n-k)\beta(k)}\right)^{q}\right)^{\frac{p}{q}%
}\right)^{\frac{1}{p}}
+\|f\|_{\beta}\|g\|_{\beta}(b_{M,1}+b_{M,2}+\cdots+b_{M,M})
\end{equation*}%
\begin{equation*}
\leq C^{\frac{1}{q}}_{o}\|g\|_{\beta}\left(\sum^{\infty}_{n=M} |\hat{f}%
(n)|^{p}\beta(n)^{p}\right)^{\frac{1}{p}}
+\|f\|_{\beta}\|g\|_{\beta}(b_{M,1}+b_{M,2}+\cdots+b_{M,M}) .
\end{equation*}%
\newline
Now, when $p=1$ we get that
\begin{equation*}
\|M_{\diamond, f}(g)-K_{M}(g)\|_{\beta}=\|f\diamond g-K_{M}(g)\|_{\beta}
\end{equation*}
\begin{equation*}
=\left
\|\sum^{\infty}_{n=M+1}\sum^{\infty}_{m=M+1}\frac{\delta_{n+m}}{%
\delta_{n} \delta_{m}}\hat{g}(n)\hat{f}(m)z^{n+m}+
\sum^{M}_{n=0}\sum^{\infty}_{m=M+1}\frac{\delta_{n+m}} {\delta_{n}\delta_{m}}%
\hat{f}(n)\hat{g}(m)z^{n+m}\right\|_{\beta}
\end{equation*}
\begin{equation*}
\leq\left(\sum^{\infty}_{n=M+1}\sum^{\infty}_{m=M+1}\frac{%
\delta_{n+m}\beta(n +m)}{\delta_{n} \delta_{m}\beta(n)\beta(m)}+
\sum^{M}_{n=0}\sum^{\infty}_{m=M+1}\frac{\delta_{n+m}\beta(n+m)} {%
\delta_{n}\delta_{m}\beta(n)\beta(m)}\right)\|f\|_{\beta}\|g\|_{\beta} .
\end{equation*}

Hence in the both case, by using (2.2) and (2.4), $\|M_{\diamond,
f}-K_{M}\|_\beta\rightarrow0$ when $M\rightarrow\infty$. Thus $M_{\diamond,
f}$ is the limit in norm of a sequence of finite-rank operators and is
therefore compact.\ $\Box$

\vspace*{0.3cm} \textbf{Lemma 2.3.} Let $f\in\ell^{p}(\beta)$ and $\hat{f}%
(0)\neq 0$. Then $M_{\diamond, f}$ is one-to-one.

\vspace*{0.3cm} \textbf{Proof.} Let $g\in\ell^{p}(\beta)$ and $M_{\diamond,
f}(g)=f\diamond g=0$. Then $\widehat{(f\diamond g)}(n)=0$, for all $n\in%
\mathbb{N}\cup\{0\}$. Since
\begin{equation*}
\widehat{({f\diamond} g)}(n)=\sum^{n}_{k=0}\frac{\delta_{n}} {%
\delta_{k}\delta_{n-k}}\hat{f}(k)\hat{g}(n-k),
\end{equation*}
we get that
\begin{equation*}
\widehat{(f\diamond g)}(0)=\hat{f}(0)\hat{g}(0)=0 \Longrightarrow\hat{g}(0)=0
\end{equation*}
\begin{equation*}
\widehat{(f\diamond g)}(1)=\hat{f}(0)\hat{g}(1)+\hat{f}(1)\hat{g}(0)=0
\Longrightarrow \hat{g}(1)=0 ,
\end{equation*}
and so on.  Thus, we get $\hat{g}(0)=\hat{g}(1)=\hat{g}(2)=\cdots =0$, and
so $g=0$.\ $\Box$

\vspace*{0.3cm} \textbf{Theorem 2.4.} If $f\in\ell^{p}(\beta)$, then f is $%
\diamond$- invertible if and only if  $\hat{f}(0)\neq 0$.

\vspace*{0.3cm} \textbf{Proof.}  Suppose that $\hat{f}(0)\neq 0$. Put $h=f-%
\hat{f}(0)$. Then $M_{\diamond, f}=\hat{f}(0)I+M_{\diamond, h}$ with $\hat{h}%
(0)=0$. By the above lemmas and the open mapping theorem, $M_{\diamond,
f}^{-1}:M_{\diamond, f}(\ell^{p}(\beta))\rightarrow \ell^{p}(\beta)$ is
bounded. On the other hand, since $M_{\diamond, f}$ is compact, then the
residual spectrum of $M_{\diamond, f}$ is empty, and so $M_{\diamond,
f}^{-1}\in{B}(\ell^{p}(\beta))$. Conversely, suppose that $f$ is invertible.
Then there exists $g\in\ell^{p}(\beta)$ such that $f\diamond g=1$ and so $%
\hat{f}(0)\hat{g}(0)=\widehat{(f\diamond g)}(0)=1$. This implies that $\hat{f%
}(0)\neq0$. \ $\Box$

\vspace*{0.3cm} \textbf{Corollary 2.5.}  The maximal ideal space of $%
(\ell^{p}(\beta),\diamond)$ consists of one homomorphism $\varphi(f)=\hat{f}%
(0)$ ($f\in\ell^{p}(\beta)$).

\vspace*{0.3cm} \textbf{Proof.} Let $\mathfrak{M}(\ell^{p}(\beta))$ be the
maximal ideal space of $\ell^{p}(\beta)$ with weighted Cauchy product $%
\diamond$. Recall that for each $f\in\ell^{p} (\beta)$, $\lambda\in \sigma(f)
$, the spectrum of $f$, if and only if $f-\lambda$ is not invertible. By
Theorem 2.4, $\lambda\in \sigma(f)$  if and only if $\widehat{(f-\lambda)}%
(0)=0$ or equivalently $\hat{f}(0)=\lambda$. Since $\sigma(f)=\{\varphi(f):%
\varphi\in \mathfrak{M}(\ell^{p}(\beta))\}$, thus $\varphi\in \mathfrak{M}%
(\ell^{p}(\beta))$ if and only if $\varphi(f)=\hat{f}(0)$, for each $%
f\in\ell^{p} (\beta)$.  \ $\Box$

\vspace*{0.3cm} Let $X$ is a Banach space and let $A\in B(X)$, the space of
all bounded linear operators on $X$. Then $x\in X$ is called cyclic vector
for $A$ in $X$ if $X=\overline{\mbox{span}}\{A^{n}x:n=0,1,2,\cdots \}$. Also
an operator $A$ in $B(X)$ is called unicellular on $X$ if the set of its
invariant closed subspaces, Lat($A$), is linearly ordered by inclusion. In
the literature there are some sufficient conditions for the usual
multiplication operator $M_{z}$ on $\ell ^{p}(\beta )$ to be unicellular
(see \cite{3,4}). In the following we give a sufficient condition for the $%
\diamond $-multiplication operator $M_{\diamond ,z}$ acting on $\ell
^{p}(\beta )$ to be unicellular.

\vspace*{0.3cm} Let $\ell^{p}_{0}(\beta)=\ell^{p}(\beta)$, $%
\ell^{p}_{\infty}(\beta)=\{0\}$ and let for $i\in\mathbb{N}\cup\{0\}$, $%
\ell^{p}_{i}(\beta)=\{\Sigma_{n\geq i} \ c_{n}z^{n}\in\ell^{p}(\beta)\}$.
Given two functions $f(z)=\Sigma^{\infty}_{n=i}  \ \hat{f}(n)z^{n}$ and $%
g(z)=\Sigma^{\infty}_{n=i}\ \hat{g}(n)z^{n}$ of the  subspace $%
\ell_{i}^{p}(\beta)$.  Let $\diamond_{i}$ be the restriction of weighted
Cauchy product $\diamond$  on $\ell_{i}^{p}(\beta)$ defined as follows:
\begin{equation*}
(f\diamond_{i}g)(z):=\sum_{n,m\geq i}\frac{\delta_{n+m-i}}{\delta_{n}
\delta_{m}}\hat{f}(n)\hat{g}(m)z^{n+m-i}.
\end{equation*}

For each $n, k, M\in\mathbb{N}\cup\{0\}$ and $1<p<\infty$, define
\begin{equation*}
C_{i}:=\sup_{n\geq i}\ \sum^{n}_{k=i}\left(\frac{\delta_{n} \beta(n)}{%
\delta_{k}\delta_{n-k+i}\beta(k)\beta(n-k+i)}\right)^{q},
\end{equation*}
and
\begin{equation*}
b^{i}_{M,k}:=\sup_{n\geq M+i+1}\frac{\delta_{n+k}\beta(n+k)} {%
\delta_{n}\delta_{k+i}\beta(n)\beta(k+i)} .
\end{equation*}
Note that when $i=0$, these are coincide with (2.1) and (2.2) respectively.
Also for $p=1$, let for some $i\leq N_{i}\in\mathbb{N}$,
\begin{equation*}
B_{i}=\sum_{n,m\geq N_{i}+1}\frac{\delta_{n+m-i}\beta(n+m-i)}{%
\delta_{n}\delta_{m}\beta(n)\beta(m)}<\infty.
\end{equation*}

For $f\in\ell_{i}^{p}(\beta)$, put

\begin{equation*}
M_{\diamond_{i}, f}(g):=f\diamond_{i}g, \ \ \ \ g\in\ell_{i}^{p}(\beta).
\end{equation*}
Then for each $n\in\mathbb{N}\cup\{0\}$ we have
\begin{equation}
M^{n}_{\diamond,z}(f)=\frac{\delta_{n}}{\delta_{1}^{n}}z^{n} \diamond f=%
\frac{\delta_{n+i}}{\delta_{1}^{n}}z^{n+i} \diamond_{i} f=M_{\diamond_{i},
f}(\frac{\delta_{n+i}}{\delta_{1}^{n}}z^{n+i}) .
\end{equation}

It follows that for each $i\in\mathbb{N}\cup\{0\}$, $\ell^{p}_{i}(\beta)$ is
an invariant subspace for $M_{\diamond,z}$. Let $\|.\|_{\beta, i}$ be the
restriction of $\|.\|_{\beta}$ on $\ell_{i}^{p}(\beta)$. Then we have
\begin{equation*}
\widehat{(f\diamond_{i}g)}(n)=\sum^{n}_{k=i}\frac{\delta_{n}}{\delta_{k}
\delta_{n-k+i}}\hat{f}(k)\hat{g}(n-k+i) ,
\end{equation*}
and $\|f\diamond_{i}g\|_{\beta,i}\leq {C_{i}}^{\frac{1}{q}%
}\|f\|_{\beta,i}\|g\|_{\beta,i}$. Thus $\ell_{i}^{p}(\beta)$ is a Banach
algebra with multiplication $\diamond_{i}$ and unity $\delta_{i}z^{i}$.
Also, by finiteness of $B_{i}$ and by the same method of proposition 2.1 $%
\ell_{i}^{1}(\beta)$ is a Banach algebra with multiplication $\diamond_{i}$
and unity $\delta_{i}z^{i}$.

\vspace*{0.3cm} \textbf{Theorem 2.6.}  Let $C_{i}<\infty $ and let $\lim
b^{i}_{M,k}=0$ when  $M\rightarrow\infty$.  Then Lat$(M_{\diamond,z})=\{%
\ell^{p}_{i}(\beta):i\geq 0\}$, and so  $M_{\diamond,z}$ is an unicellular
operator on  $(\ell^{p}(\beta), \diamond)$.

\vspace*{0.3cm} \textbf{Proof.} Note that the set $\{\ell^{p}_{i}(\beta):i%
\geq 0\}$ is linearly ordered by inclusion.

For $f\in\ell^{p}(\beta)$, put
\begin{equation*}
E(f)=\overline{\mbox{span}}\{f, \ M_{\diamond,z}(f), \
M_{\diamond,z}^{2}(f), \ M_{\diamond,z}^{3}(f), \ldots\} .
\end{equation*}
Since $E(f)$ is an invariant subspace for $M_{\diamond,z}$, the operator $%
M_{\diamond,z}$ is unicellular in the Banach algebra $\ell^{p}(\beta)$ if
and only if for all nonzero $f\in\ell^{p}(\beta)$, $E(f)$ is equal to $%
\ell_{i}^{p}(\beta)$ for some $i= 0, 1, 2, \ldots, i(f)$.

\vspace*{0.3cm} Now, we show that equality of $E(f)$ with $%
\ell_{i}^{p}(\beta)$ is equivalent to the condition $\hat{f}(i)\neq 0$. To
see this note that by (2.5) we have
\begin{equation*}
E(f)=\overline{\mbox{span}}\{M^{n}_{\diamond,z}(f) : n\geq 0\}= \overline{%
\mbox{span}} \{M_{\diamond_{i}, f}(\frac{\delta_{n+i}}{\delta_{1}^{n}}%
z^{n+i}) : n\geq 0\} = \overline{M_{\diamond_{i}, f}(\ell^{p}_{i}(\beta))} ,
\end{equation*}
and so
\begin{equation*}
E(f)=\ell^{p}_{i}(\beta) \Longleftrightarrow \overline{M_{\diamond_{i},
f}(\ell^{p}_{i} (\beta))}=\ell^{p}_{i}(\beta).
\end{equation*}

Now, we show that
\begin{equation*}
\overline{M_{\diamond_{i}, f}(\ell^{p}_{i}(\beta))}=\ell^{p}_{i}(\beta)%
\Longleftrightarrow \hat{f}(i)\neq0.
\end{equation*}
Indeed, if $\overline {M_{\diamond_{i}, f}(\ell^{p}_{i}(\beta))}%
=\ell^{p}_{i}(\beta)$ then there exists a sequence of $\{f_{n}\}\subseteq%
\ell^{p}_{i}(\beta)$  such that $f\diamond_{i}f_{n}\rightarrow\delta_{i}z^{i}
$, and so $\frac{1}{\delta_{i}}\hat{f}(i)\hat{f_{n}}(i)\rightarrow \delta_{i}
$ as $n\rightarrow\infty$,  which implies that $\hat{f}(i)\neq0$.

\vspace*{0.3cm} Conversely, if $\hat{f}(i)\neq0$ then it is sufficient prove
that $M_{\diamond_{i}, f}$ is an invertible operator in  $\ell_{i}^{p}(\beta)
$  which will imply that $\overline{M_{\diamond_{i}, f}(\ell^{p}_{i}(\beta))}%
=\ell^{p}_{i}(\beta)$,  as desired. Put $h=f-\hat{f}(i)z^{i}$. Then $%
M_{\diamond_{i}, f}=\frac{\hat{f}(i)}{\delta_{i}}I+M_{\diamond_{i}, h} $. By
the same argument in the proof of the Lemmas 2.2 and 2.3, it is sufficient
to show that $M_{\diamond_{i}, h}$ and $M_{\diamond_{i}, f}$ are compact and
one-to-one  operators respectively in $\ell_{i}^{p}(\beta)$. For this
purpose define
\begin{equation*}
K^{i}_{M}(g)=\sum^{i+M}_{m=i}\sum^{i+M}_{n=i} \frac{\delta_{m+n-i}}{%
\delta_{n}\delta_{m}} \hat{h}(n)\hat{g}(m)z^{n+m-i}+ \sum^{i+M}_{m=i}\frac{%
\hat{g}(m)\delta^{m-i}_{1}}{\delta_{m}}R_{M+1}(M^{m-i}_{\diamond,z}(h)),
\end{equation*}
for every integer $M\in\mathbb{N}$ and $g\in\ell_{i}^{p}(\beta)$. By an
argument similar to the proof of the Lemma 2.2, $K^{i}_{M}$ is a finite-rank
operator and
\begin{equation*}
\|M_{\diamond_{i}, h}(g)-K^{i}_{M}(g)\|_{\beta,i}=
\end{equation*}%
\begin{equation*}
\leq{C_{i}}^{\frac{1}{q}}\|g\|_{\beta,i}(\sum^{\infty}_{n=M+i+1}| \hat{h}%
(n)|^{p}\beta(n)^{p})^{\frac{1}{p}}
+\|h\|_{\beta,i}\|g\|_{\beta,i}(b^{i}_{M,1}+b^{i}_{M,2}+\cdots +b^{i}_{M,M}).
\end{equation*}

Also, when $p=1$ we get that
\begin{equation*}
\|M_{\diamond_{i}, h}(g)-K^{i}_{M}(g)\|_{\beta,i}
\end{equation*}
\begin{equation*}
\leq\left\{\left(\sum^{\infty}_{n=i+M+1}\sum^{\infty}_{m=i+M+1}+
\sum^{i+M}_{n=i}\sum^{\infty}_{m=i+M+1} \right) \frac{\delta_{n+m-i}\beta(n
+m-i)}{\delta_{n} \delta_{m}\beta(n)\beta(m)}\right\}\|h\|_{\beta,i}\|g\|_{%
\beta,i} .
\end{equation*}

Thus for $1\leq p<\infty$, $\|M_{\diamond_{i}, h}-K^{i}_{M}\|\rightarrow0$
when $M\rightarrow\infty$, and hence $M_{\diamond_{i}, h}$ is a compact
operator. Now, let $g\in\ell_{i}^{p}(\beta)$ and $M_{\diamond_{i},
f}(g)=f\diamond_{i} g=0$. Then for all $n\geq i$, $\widehat{(f\diamond_{i} g)%
}(n)=0$. Hence we get that
\begin{equation*}
\widehat{(f\diamond_{i} g)}(i)=\frac{1}{\delta_{i}}\hat{f}(i)\hat{g}(i)=0
\Longrightarrow \hat{g}(i)=0
\end{equation*}
\begin{equation*}
\widehat{(f\diamond_{i} g)}(i+1)=\frac{1}{\delta_{i}}\hat{f}(i) \hat{g}(i+1)+%
\frac{1}{\delta_{i}}\hat{f}(i+1)\hat{g}(i)=0 \Rightarrow \hat{g}(i+1)=0 ,
\end{equation*}
and so on. It follows that $\hat{g}(i)=\hat{g}(i+1)=\hat{g}(i+2)=\cdots = 0$%
; i.e., $g=0$. This completes the proof of the theorem. \ $\Box$

\vspace*{0.3cm} \textbf{Corollary 2.7.} Let $f\in\ell^{p}(\beta)$. Then $f$
is a cyclic vector for $M_{\diamond,z}$ if and only if $\hat{f}(0)\neq0$.

\vspace*{0.3cm} \textbf{Proof.} Let $f\in\ell^{p}(\beta)$. Then we have
\begin{equation*}
\overline{\mbox{span}}\{M^{n}_{\diamond, z}(f) : n\geq 0 \}=\overline{%
\mbox{span}} \{M_{\diamond, f}(\frac{\delta_{n}}{\delta^{n}_{1}}z^{n}) :
n\geq 0 \}.
\end{equation*}
If $\hat{f}(0)\neq0$, then by Theorem 2.5, $M_{\diamond, f}$ is an
invertible operator on $\ell^{p}(\beta)$, and so $\overline{\mbox{span}}
\{M^{n}_{\diamond, z}(f) : n\geq 0 \}= \ell^{p}  (\beta)$ which implies that
$f$ is a cyclic vector for $M_{\diamond, z}$.  Conversely, suppose $f$ is a
cyclic vector for $M_{\diamond, z}$.  Then there exists sequence $%
\{f_{n}\}\subseteq\ell^{p}(\beta)$  such that $\|f_{n}\diamond
f-1\|_\beta\rightarrow 0$. This implies that  $\hat{f_{n}}(0)\hat{f}%
(0)\rightarrow1$, and so $\hat{f}(0)\neq 0$. \ $\Box$

\vspace*{0.3cm} In the following theorem we characterize the form of all
closed ideals of the Banach algebra $(\ell^{p}(\beta),\diamond)$.

\vspace*{0.3cm} \textbf{Theorem 2.8.} Let $i\in\mathbb{N}\cup\{0\}$, $%
C_{i}<\infty $ and let $\lim b^{i}_{M,k}=0$ when  $M\rightarrow\infty$. Then
the closed ideals of $(\ell^{p}(\beta),\diamond)$ are exactly of the form $%
\ell^{p}_{i}(\beta)$.

\vspace*{0.3cm} \textbf{Proof.} For any $i\in\mathbb{N}\cup\{0\}$, it is
easy to see that $\ell^{p}_{i}(\beta)$ is a closed ideal of $%
(\ell^{p}(\beta),\diamond)$. Let $K$ be an arbitrary closed ideal of $%
(\ell^{p}(\beta),\diamond)$. Then for each $f\in K$, $z\diamond f\in K$ and
so $M_{\diamond, z}(K)\subseteq K$. Now, by Theorem 2.6, $%
K=\ell^{p}_{i}(\beta)$ for some $i\in\mathbb{N}\cup\{0\}$. This completes
the proof. \ $\Box$\

\end{document}